\documentclass[12pt]{amsart}

\newtheorem{theorem}{Theorem}[section]
\newtheorem{lemma}[theorem]{Lemma}

\theoremstyle{definition}   
\newtheorem{definition}{Definition}

\theoremstyle{remark}

\numberwithin{equation}{section}

\usepackage{times}
\usepackage{enumerate}
\usepackage{tfrupee}
\usepackage{tikz}
\usetikzlibrary{chains,fit}
\usetikzlibrary{shapes,snakes}
\usetikzlibrary{graphs}
\usepackage{graphicx,adjustbox}
\usepackage{hyperref}
\usepackage{amsmath,amssymb}

\newcommand{\LT}{\mbox{\rm Lt}}

\usepackage{amscd}
\usepackage{graphicx}
\usepackage[all]{xy}

\title[ASL Structures]
{ASL structures of some quadrics}
\author{
Joydip Saha
\and
Indranath Sengupta
}
\date{}

\address{\small \rm  Stat-Math Unit, Indian Statistical Institute, 203 B.T. Road, Kolkata 700 108.} 
\email{saha.joydip56@gmail.com}
\thanks{The first author is supported by the NPDF fellowship PDF/2019/001074, 
sponsored by the SERB, Government of India.}

\address{\small \rm  Discipline of Mathematics, IIT Gandhinagar, Palaj, Gandhinagar, 
Gujarat 382355, INDIA.}
\email{indranathsg@iitgn.ac.in}
\thanks{The second author is the corresponding author. This research is supported by the 
MATRICS research grant MTR/2018/000420, sponsored by SERB.}
\date{}



\begin{document}
\begin{abstract}
Let $K$ be a field and $X$, $Y$ denote matrices such that, the entries 
of $X$ are either indeterminates over $K$ or $0$ and the entries 
of $Y$ are indeterminates over $K$ which are different from those appearing 
in $X$. We consider ideals of the form $I_{1}(XY)$, which is the 
ideal generated by the $1\times 1$ minors of the matrix $XY$. We prove 
that the quotient ring $K[X, Y]/I_{1}(XY)$ admits an ASL structure for 
certain $X$ and $Y$.
\end{abstract}

\maketitle

\section*{Introduction}

Let $K$ be a field and $\{x_{ij}; \, 1\leq i \leq m, \, 1\leq j \leq n\}$, 
$\{y_{j}; \, 1\leq j \leq n\}$ be indeterminates over $K$. Let $R = K[x_{ij}]$ 
and $S = K[x_{ij}, y_{j}]$ denote the polynomial algebras over 
$K$. Let $X$ denote an $m\times n$ matrix such that its entries belong to the 
ideal $\langle \{x_{ij}; \, 1\leq i \leq m, \, 1\leq j \leq n\}\rangle$ and 
$Y=(y_{j})_{n\times 1}$ the generic $n\times 1$ column matrix. Let 
$I_{1}(XY)$ denote the ideal generated by the $1\times 1$ minors or the 
entries of the $m\times 1$ matrix $XY$. We assume that $m=n$ and write 
$\mathcal{I} =  I_{1}(XY) = \langle g_{1}, g_{2}, \ldots , g_{n} \rangle$. 
The ideal $I_{1}(XY)$ is a special case of the defining ideal of a variety of
complexes, see \cite{concini}. These ideals also feature in \cite{tchernev}, 
in the study of the structure of a \textit{universal ring} of a 
\textit{universal pair}. Tchernev has proved that the set of standard 
monomials form a free basis for the universal ring. 
\medskip

Under the assumption that $X$ is genereic (respectively generic symmetric) 
and with respect to any monomial order satisfying $x_{11}> x_{22}> \cdots >x_{nn}$; 
$x_{ij}, y_{j} < x_{nn}$ for every $1 \leq i \neq j \leq n$ (respectively 
$1 \leq i < j \leq n$), it is true that the set $\{g_{1}, g_{2}, \ldots , g_{n}\}$ 
forms a Gr\"{o}bner basis for the ideal $I_{1}(XY)$; see \cite{sst}. Another 
Gr\"{o}bner basis exists and that also appears in \cite{sst}, 
which has been used to prove normality in \cite{sst} and compute primary decomposition 
of these ideals in \cite{sstp}. In this paper we will show that the knowledge of 
Gr\"{o}bner basis actually leads us to the fact that $K[X, Y]/I_{1}(XY)$ admits an 
ASL structure.   

\section{Algebra with Straightening Law}
ASL or the \textit{Algebra with Starightening Law} is a special structure on an 
algebra $A$, over a partially ordered subset of $A$. We list the definition and 
some basic facts below. We refer to \cite{eisenbud} for definitions and pertinent results.
 
\begin{definition}
Let $A$ be a commutative ring with $1$ and $H$ a subset of $A$. Suppose that $H$ 
is a partially ordered set (poset). A \textit{standard monomial} is a product of 
the form $m_{1}\cdots m_{k}$, such that $m_{1}\leq\cdots\leq m_{k}$.  
\end{definition}

\begin{definition}
Let $B$ be a commutative ring with $1$ and  $A$ an algebra over the ring $B$. 
Let $H$ be a finite partially ordered subset of $A$, which generates $A$ as 
a $B$ algebra. Then $A$ is an algebra with straightening law (on $H$, over $B$) 
if the following conditions are satisfied.
\begin{enumerate}
\item The algebra $A$ is a free $B$ module whose basis is the set of standard monomials.
\item If $\alpha$ and $\beta$ in $H$ are incomparable and if 
$$\alpha\beta=\sum r_{i}m_{i_{1}}m_{i_{2}}\cdots m_{i_{k_{i}}},$$ 
where $r_{i}\neq 0$ is in $B$ and $m_{i_{1}}\leq m_{i_{2}}\leq \cdots \leq m_{i_{k_{i}}}$ 
is the unique expression for $\alpha\beta$ in $A$ as a linear combination of standard monomials, 
then $m_{i_{1}}\leq \alpha,\beta$ for every $i$. 
\end{enumerate} 
\end{definition}

\begin{theorem}
Let us consider the $K$ algebra $S/I_{1}(XY)$, where $X=(x_{ij})$ is generic $n\times n$ matrix of 
indeterminates $x_{ij}$ and $Y$ is generic $n\times 1$ matrix of indeterminates $y_{i}$. Then 
$S/I_{1}(XY)$ is an algebra with straightening law on the partially ordered set 
$H=\{x_{ij}+I_{1}(XY),y_{i}+I_{1}(XY)\mid 1\leq i,j\leq n \}$ over $K$. The partial order 
$\preceq$ on $H$ is given by following chains: 
\begin{enumerate}
\item $\bar{x}_{12}\preceq\ldots \preceq\bar{x}_{1n}\preceq\bar{x}_{21}\preceq\bar{x}_{23}\preceq\ldots \preceq\bar{x}_{2, n}\preceq\ldots\preceq\bar{x}_{n, (n-1)}\preceq \bar{x}_{nn}\preceq \bar{x}_{(n-1), (n-1)}\preceq \ldots \bar{x}_{11}$,
\item $\bar{x}_{n, (n-1)}\preceq\bar{y}_{n}\preceq\ldots\preceq\bar{y}_{1}$,
\item $\bar{x}_{22}\preceq\bar{y}_{1}$,
\item $\bar{y}_{n}\preceq \bar{x}_{(n-1), (n-1)} $,
\item $\bar{x}_{(i+1), (i+1)}\bar{y}_{i}\preceq \bar{x}_{(i-1), (i-1)}$, for $2\leq i\leq n-1$.\\
\end{enumerate}
Here $\bar{h}$ denotes the residue modulo $I_{1}(XY)$.

\end{theorem}

\proof We fix monomial order on $S$ as in the theorem for the Gr\"{o}bner basis 
for the generic case:

\begin{enumerate}
\item $x_{11}> x_{22}> \cdots >x_{nn}$;
\item $x_{ij}, y_{j} < x_{nn}$ for every $1 \leq i \neq j \leq n$.
\end{enumerate}

\noindent Let $\mathcal{I} = I_{1}(XY) = \langle g_{1},g_{2},\cdots ,g_{n}\rangle$, 
where $g_{i}=\sum_{j=1}^{n} x_{ij}y_{j}$. Then, the set $\lbrace g_{1},\cdots ,g_{n}\rbrace$ 
forms a Gr\"{o}bner basis for the ideal $\mathcal{I}$ with respect to the monomial order 
written above. Therefore, $\textrm{in}(\mathcal{I})=\langle\{x_{ii}y_{i}\mid 1\leq i\leq n\}\rangle$ and the set $L=\{\bar{m}\mid m \,\textrm{is a monomial such that } m\notin \textrm{in}(\mathcal{I})\}$ 
forms a basis of $K$ algebra $S/\mathcal{I}$. Since only $\bar{x}_{ii}$ and $\bar{y}_{i}$ are incomparable 
in $H$, for all $1\leq i\leq n$, then it is obvious that $L$ is set of standard monomials in $S/\mathcal{I}$ 
with respect to the given partial order on $H$. Therefore the first condition in ASL holds. 
Now we have the expression $$\bar{x}_{ii}\bar{y}_{i}=-(\sum_{j=1,\,j\neq i}^{n} \bar{x}_{ij}\bar{y}_{j}).$$ 
Here, for each $1\leq i\leq n$, we have $\bar{x}_{ij}\preceq \bar{y}_{j} $ and 
$\bar{x}_{ij}\preceq \bar{y}_{i}$, $\bar{x}_{ij}\preceq \bar{x}_{ii}$, for all $j\neq i$.\qed

\bibliographystyle{amsalpha}

\end{document}